\documentclass{amsproc}
\usepackage[all]{xy}

\newtheorem{theorem}{Theorem}[section]
\newtheorem{proposition}[theorem]{Proposition}

\theoremstyle{definition}
\newtheorem{definition}[theorem]{Definition}
\newtheorem{example}[theorem]{Example}

\theoremstyle{remark}
\newtheorem{remark}[theorem]{Remark}

\numberwithin{equation}{section}

\newcommand{\comment}[1]{}

\newcommand{\Rr}{\mathbb R}
\newcommand{\Zz}{\mathbb Z}

\newcommand{\set}[1]{\left\{#1\right\}}

\newcommand{\To}{\longrightarrow}

\newcommand{\X}{\ensuremath{\mathcal{X}}}
\newcommand{\F}{\ensuremath{\mathcal{F}}}

\newcommand{\Agerm}{\mathfrak{Aut}}
\newcommand{\Ogerm}{\mathfrak{Out}}

\newcommand{\Aut}{\text{\rm Aut}\,}
\newcommand{\Inn}{\text{\rm Inn}\,}
\newcommand{\End}{\text{\rm End}\,}
\newcommand{\Vect}{\text{\rm Vect}}
\newcommand{\Vectd}{\text{\rm Vect}^{\text{def}}}
\newcommand{\gl}{\mathfrak{gl}}

\newcommand{\G}{\mathcal{G}}            

\newcommand{\A}{A}                      
\newcommand{\al}{\alpha}                
\newcommand{\Lie}{\mathcal{L}}          
\renewcommand{\gg}{\mathfrak{g}}        
\newcommand{\Ker}{\text{\rm Ker}\,}     
\renewcommand{\Im}{\text{\rm Im}\,}     
\newcommand{\Ad}{\text{\rm Ad}\,}       
\newcommand{\ad}{\text{\rm ad}\,}       

\newcommand{\rank}{\text{\rm rk}\,}   

\begin{document}

\title{Invariants of Lie algebroids}

\author{Rui Loja Fernandes}
\address{Depart.~de Matem\'{a}tica, Instituto Superior T\'{e}cnico,
1049-001 Lisboa, PORTUGAL}
\email{rfern@math.ist.utl.pt}
\thanks{Supported in part by FCT through program POCTI and grant
  POCTI/1999/MAT/33081.} 

\date{October, 2001}
\begin{abstract}
Several new invariants of Lie algebroids have been discovered
recently. We give an overview of these invariants and 
establish several relationships between them.
\end{abstract}
\maketitle

\section{Introduction}             %
\label{Invariants:Outline}         %
%
It is becoming increasingly apparent that Lie algebroids provide the
appropriate setting for developing the differential geometry of singular
geometric structures. The study of global properties of Lie algebroids
is therefore a way of approaching the global theory of singular
geometric structures, about which little is known.  In this survey we
shall describe several constructions of Lie algebroid invariants which
have been introduced in the last few years. Although we are just
starting to grasp their properties, it is clear that they play an
important role in understanding the global behavior of singular
geometric structures.

Let $\pi:A\to M$ be a Lie algebroid over $M$ with anchor $\#:A\to TM$
and Lie bracket $[~,~]:\Gamma(A)\times\Gamma(A)\to\Gamma(A)$. For
general definitions and conventions we refer the reader to
\cite{CaWe}, and further background material is given in
\cite{Mack1,Mack2}. We shall adopt the point of view that such a Lie
algebroid describes a concrete geometric situation. In each case, the
Lie algebroid plays the role of the \emph{tangent bundle}, and many
constructions to be given later can be traced back to this simple
idea. Some typical cases we have in mind are:
\vskip 10 pt

\subsubsection*{Ordinary Geometry:} We take $A=TM$, $\#$ the
identity, and $[~,~]$ the usual Lie bracket of vector fields. Later,
we will see that virtually every construction for an arbitrary Lie
algebroid reduces to some well-known construction when applied to this
case.
\vskip 10 pt

\subsubsection*{Lie Theory:} At the other extreme, we take $A=\gg$ to be a
Lie algebra and $M=\{*\}$ a single point set. Although this is a
somewhat degenerate case, it is useful because the local-global
dichotomy for Lie algebroids resembles Lie theory. Also, the
terminology is usually motivated either by this case or the
previous one.
\vskip 10 pt

\subsubsection*{Equivariant Geometry:} If we are given an action of a
Lie algebra on a manifold, i.e., a homomorphism $\rho:\gg\to \X(M)$ to
the Lie algebra of vector fields on M, then we have a naturally
associated action Lie algebroid. The bundle $A$ is the trivial vector
bundle $M\times\gg\to M$, the anchor $\#:M\times\gg\to TM$ is defined
by
\[ \#(x,v)\equiv\rho(v)|_x,\]
and the Lie bracket is given by
\[ [v,w](x)=[v(x),w(x)]+(\rho(v(x))\cdot w)|_{x}-(\rho(w(x))\cdot v)|_{x},\]
where we identify a section $v$ of $M\times\gg\to M$ with a
$\gg$-valued function $v:M\to\gg$. In this case we have a particularly
simple geometric interpretation for the orbit foliation of the
algebroid\footnote{In this paper, foliations can be \emph{singular}
as in Sussmann \cite{Suss}. By a \emph{regular foliation} we mean a
non-singular foliation.}, even though the action does not 
always integrate to a global Lie group action.
\vskip 10 pt

\subsubsection*{Foliation Theory:} Let $\F$ be a regular foliation of
$M$. The associated involutive distribution $A=T\F$ has a Lie
algebroid structure with anchor the inclusion into $TM$ and bracket
the Lie bracket of tangent vector fields to $\F$. Many constructions
in Lie algebroid theory, related to the geometry and topology of the
orbit foliation, are inspired by constructions in foliation theory.
\vskip 10 pt

\subsubsection*{Poisson Geometry:} Consider a Poisson manifold
$(M,\pi)$, where $\pi\in\Gamma(\bigwedge^2 TM)$ is a bi-vector field
satisfying $[\pi,\pi]=0$. It is well known that the cotangent bundle 
$A=T^*M$ has a natural Lie algebroid structure, where the anchor
$\#:T^*M\to TM$ is contraction with $\pi$, and the bracket on 1-forms 
is the Koszul bracket:
\[ [\al,\beta]=\Lie_{\#\al} \beta-\Lie_{\#\beta}\al-d\pi(\al,\beta).\]
Many concepts to be discussed below were first introduced for Poisson 
manifolds and then generalized to Lie algebroids.
\vskip 20 pt

As a first example of an invariant let us consider Lie algebroid
cohomology (see \cite{Mack1} for more details). We just mimic the
usual definition of de Rham cohomology: the space of differential forms is
$\Omega^r(A)=\Gamma(\wedge^r A^*)$, and we define the exterior
differential $d_A:\Omega^\bullet(A)\to\Omega^{\bullet+1}(A)$ by:
\begin{multline}
\label{eq:differential}
d_A Q(\al_0,\dots,\al_r)=\frac{1}{r+1}\sum_{k=0}^{r+1}
(-1)^k\#\al_k(Q(\al_0,\dots,\widehat{\al}_k,\dots,\al_r)\\
+\frac{1}{r+1}\sum_{k<l}(-1)^{k+l+1}Q([\al_k,\al_l],\al_0,\dots,\widehat{\al}_k,\dots,\widehat{\al}_l,\dots,\al_r).
\end{multline}
where $\al_0,\dots,\al_r$ are any sections of $A$. In this way we
obtain a complex $(\Omega^\bullet(A),d_A)$, and the corresponding
cohomology is called the \emph{Lie algebroid cohomology} of $A$ 
(with trivial coefficients) and
denoted $H^\bullet(A)$. The (dual of the) anchor induces a map 
\[\#^*:H^\bullet_{\text{de Rham}}(M)\to H^\bullet(A),\]
which is usually neither injective nor surjective. For the geometric
situations discussed above we obtain well-known cohomology theories
such as de Rham cohomology, Lie algebra cohomology, foliated
cohomology and Poisson cohomology.  As these examples show, Lie
algebroid cohomology may not be homotopy invariant and hence it may
be hard to compute (to say the least).

The problem of computing Lie algebroid cohomology is intimately
related with the singular behavior of the orbit foliation of the Lie
algebroid. The same will be true about all other invariants to be
introduced below, and this is in fact one of the main topics of the
present work.  The construction of the new invariants resembles the
construction of Lie algebroid cohomology, in as much as, if one knows
the proper conceptual general definitions, then the appropriate
construction will be similar to the corresponding construction in
standard geometry.

The plan of this paper is as follows.  First we consider the
fundamental group(oid) of a Lie algebroid (Section 2), which was
introduced in \cite{CrFe} for the purpose of integrating Lie
algebroids to Lie groupoids, and which is inspired by the construction
of the fundamental group(oid) of a manifold. Then we consider non-linear
and linear holonomy (Sections 3 and 4) in the spirit of foliation
theory, which was defined in \cite{Fer1} for Lie algebroids. These
will lead us naturally to primary and secondary characteristic classes
(Section 5) for Lie algebroids, which were introduced in
\cite{Cra,Fer1,Ku}. The last invariant we shall discuss is K-theory
(Section 6) which was introduced in \cite{Ginz}, and may be considered
as an extension of ordinary topological $K$-theory.

\begin{center}
\textit{Acknowledgments}
\end{center}

Parts of this manuscript formed the subject of a micro-course that I gave
at the 2001 Krynica meeting in Geometry and Topology of Manifolds. I
would like to thank Jan Kubarski, Tomasz Rybicki and Robert Wolak for
the invitation, which provided an incentive to writing this survey.

Much of the work reported here is joint work with Marius Crainic, to
whom I am in debt, and who also made significant contributions on his
own. Many other people have also influenced my view on this subject,
but I would especially like to thank Ana Cannas da Silva, Viktor
Ginzburg, Kiril Mackenzie and Alan Wenstein, for many enlightening
discussions, and Roger Picken and the anonymous referee for pointing
out many typos and suggesting many improvements to the first version
of the manuscript.

\section{The Weinstein groupoid}                           %
\label{Weinstein:groupoid}                                 %

In \cite{CrFe}, for \emph{any} Lie algebroid $A$ we have constructed
a topological groupoid $\G(A)$, called the \emph{Weinstein groupoid}, and
which is a fundamental invariant of $A$. It should be
thought of as the ``monodromy groupoid'' or ``fundamental groupoid"
of $A$: it is given as the set of equivalence classes of
$A$-paths under $A$-homotopy
\begin{equation}
\label{eq:weinstein:groupoid}
\G(A)=P(A)/\sim\ ,
\end{equation}
where:
\begin{itemize}
\item $P(A)$ denotes the set of \emph{$A$-paths}, i.e., paths 
$a:I\to A$ on the interval $I=[0,1]$, such that
\[ \frac{d}{dt}\pi(a(t))=\#a(t). \]
An $A$-path can be identified with a Lie algebroid morphism $TI\to
A$. From a differential-geometric point of view, these are
precisely the paths along which parallel transport can be performed
whenever a connection has been chosen (see below and
\cite{Fer1,Fer2}).
\item $\sim$ denotes an equivalence relation, called
\emph{homotopy of $A$-paths}, that can be
described at an abstract level as follows: two $A$-paths $a_0$
and $a_1$ are homotopic iff there exists a Lie algebroid
homomorphism $T(I\times I)\to A$, which covers a (standard)
homotopy between the base paths $\pi(a_i(t))$, with fixed end-points,
and which restricts to $a_i(t)$ on the boundaries.
\end{itemize}
The groupoid structure on $\G(A)$ comes from concatenation of
$A$-paths, which becomes associative when one passes to the
quotient. The structure maps are the obvious ones (as in the
monodromy groupoid).

To be able to work with $\sim$ one needs more concrete descriptions
which are furnished in \cite{CrFe}. There we show, for example, that this
equivalence relation is the orbit equivalence relation of a Lie
algebra action. Namely, the Lie algebra of time-dependent sections
of $A$, vanishing at the end-points,
\[ P_{0}\Gamma(A)= \set{ I\ni t\mapsto \eta_{t}\in \Gamma(A): \eta_{0}
                 = \eta_{1}=0, \eta\text{ is of class } C^2\text{ in }
t} \] 
acts on $P(A)$. So we have a Lie algebra homomorphism
\[  P_{0}\Gamma(A) \to \X(P(A)),\quad \eta\mapsto X_{\eta} \]
for which the image is precisely the tangent space to the orbits of $\sim$.

Clearly, one cannot expect the differentiable structure on the path
space to go over to the quotient. One can show (see \cite{CaFe,CrFe})
that this action gives a smooth foliation on the Banach manifold
$P(A)$, for which the orbits are smooth submanifolds of finite codimension
equal to $\dim M+\rank A$. So the most one can say, in general, is that the
Weinstein groupoid is of the same topological type as the orbit space
of a foliation.

In order to give the precise obstructions for the Weinstein groupoid
to be a Lie groupoid, and hence also the obstructions to integrating a
Lie algebroid, we introduce certain \emph{monodromy groups} of the Lie
algebroid. For that purpose, observe that an element in the isotropy
Lie algebra $\gg_x= Ker(\#_{x})$ determines a constant $A$-path, and
so we can set:

\begin{definition}
For each $x\in M$, the \textsc{monodromy group} based at $x$ is the
subgroup $N_{x}(\A)\subset \A_x$ consisting of those elements $v\in
Z(\gg_x)$ which are homotopic to zero as $A$-paths.
\end{definition}

Since the monodromy group $N_x(A)$ lies in the center $Z(\gg_x)$ of
the isotropy Lie algebra, we can identify it with an abelian subgroup
of the simply connected Lie group $\G(\gg_x)$ integrating the Lie
algebra $\gg_x$. Henceforth, we use this identification with no
further comment.

The obstructions to integrability are related to the lack of
discreteness of the monodromy groups. To explain this, let us observe
that the monodromy based at $x$ arises as the image of a second order
monodromy map
\[ \partial: \pi_{2}(L,x)\to \G(\gg_x)\]
which relates the topology of the leaf $L$ through $x$ with the 
simply-connected Lie group $\G(\gg_x)$ integrating the isotropy Lie algebra
$\gg_x= Ker(\#_{x})$. From a conceptual point of view, the monodromy
map can be viewed as an analogue of a boundary map of the homotopy
long exact sequence of a fibration. Namely, if we consider the short
exact sequence
\[ 0\To \gg_{L}\To A_{L}\stackrel{\#}{\To} TL\To 0\] 
as analogous to a fibration, the first few terms of the associated long
exact sequence will be
\[ \cdots\To\pi_{2}(L,x)\stackrel{\partial}{\To} \G(\gg_x)
\To \G(A)_x\To \pi_1(L,x) .\]
We have shown in \cite{CrFe} that $\Im \partial$ lies in the center of
$\G(\gg_x)$ and its intersection with the connected component of the
identity of $Z(\G(\gg_x))$ coincides with $N_x(A)$. With these
notations we have the following fundamental result:

\begin{theorem}[\textsc{Obstructions to Integrability \cite{CrFe}}]
\label{thm:integrability}
For a Lie algebroid $\A$ over $M$, the fol\-low\-ing are equivalent:
\begin{enumerate}
\item[(i)] $\A$ is integrable;
\item[(ii)] The monodromy groups are uniformly discrete.
\end{enumerate}
\end{theorem}

Let us be more precise about (ii). In order to measure the
discreteness of the groups $N_{x}(\A)$, set 
\[ r(x)= d(0, N_{x}(\A)-\{0\})\]
where the distance is computed with respect to an arbitrary norm on 
the vector bundle $\A$ and we adopt the convention
$d(0,\emptyset)=+\infty$. Notice that $N_{x}(A)\subset \A_x$ is
discrete iff $r(x)>0$. Then condition (ii) can be stated as
\begin{enumerate}
\item[(iia)] For all $x\in M$, $r(x)>0$;
\item[(iib)] For all $x\in M$, $\liminf_{y\to x}r(y)>0$.
\end{enumerate}
Since the monodromy groups $N_x(A)$ are isomorphic as $x$ varies in a
leaf, (iia) is an obstruction along the leaves, while (iib) is 
an obstruction transverse to the leaves.

These obstructions are computable in many examples. Given any 
splitting $\sigma: TL\to A_L$ of the short
exact sequence above, the \emph{curvature} of $\sigma$ is the
$\gg_L$-valued 2-form $\Omega\in \Omega^{2}(L;\gg_L)$ defined by:
\[ \Omega(X, Y)\equiv \sigma ([X, Y])- [\sigma(X), \sigma(Y)]\ .\]
Assume that there exists a splitting such that this 2-form takes 
values in the center $Z(\gg_L)$. Then the monodromy map 
$\partial:\pi_{2}(L,x)\to\nu^*(L)$ is given by
\begin{equation}
\label{eq:monodromy}
\partial([\gamma])=\int_{\gamma} \Omega,
\end{equation}
and this gives an effective procedure to compute the monodromy in many
example (see \cite{CrFe}, section 3.4).  Note that in this case,
$Z(\gg_{L})$ is canonically a flat vector bundle over $L$. The
corresponding flat connection can be expressed with the help of the
splitting $\sigma$ as
\begin{equation}
\label{eq:bott}
\nabla_{X}\alpha= [\sigma(X), \alpha ],
\end{equation}
and it is easy to see that the definition does not depend on
$\sigma$. In this way $\Omega_{\sigma}$ appears as a $2$-cohomology
class with coefficients in the local system defined by $Z(\gg_{L})$
over $L$, and then the integration (\ref{eq:monodromy}) is just the
usual pairing between cohomology and homotopy. In practice one can
always avoid working with local coefficients: if $Z(\gg_{L})$ is not
already trivial as a vector bundle, one can achieve this by pulling
back to the universal cover of $L$ (where parallel transport with
respect to the flat connection gives the desired trivialization).

The connection given by (\ref{eq:bott}) is the \emph{Bott
connection} on $A_L$, which can be introduced even for singular
leaves. This connection will be discussed further below, when we
consider the theory of holonomy for Lie algebroids, therefore
providing a relation between these two invariants. 

Many, if not all, results on integrability of Lie algebroids (see, e.g.,
\cite{AlHe,DaHe,DoLa,MaXu,Ni,Wein}) follow from Theorem
\ref{thm:integrability}. We give a few examples and refer the reader
to \cite{CrFe} for further examples and details.

\begin{example}
Let $A=TM$ be the tangent bundle Lie algebroid structure. A $TM$-path
is just an ordinary path in $M$ (given by its derivative), and a
$TM$-homotopy is an ordinary homotopy in $M$ with fixed end points, so
we have
\[ 
\G(TM)=\set{(x,[\gamma],y): x,y\in M, \gamma\text{ is a path from
    $x$ to $y$}}.
\] 
Therefore, $\G(TM)$ is just the fundamental groupoid $\pi_1(M)$.

More generally, for any foliation $\F$ of $M$, we have the fundamental
groupoid $\pi_1(\F)$, where now $T\F$-paths are just $\F$-paths (paths
lying on any fixed leaf) and $T\F$-homotopies are homotopies within the set
of $\F$-paths.

In both these cases the obstructions obviously vanish, and this
corresponds to the well-known fact that $\pi_1(M)$ and $\pi_1(\F)$ are
Lie groupoids.
\end{example}

\begin{example}
Let $\gg$ be a Lie algebra. Again, there are no obstructions and
Theorem \ref{thm:integrability} gives Lie's third theorem. The
construction of $\G(\gg)$ given above coincides with the construction
of the simply-connected Lie group integrating $\gg$ which is given in the
recent monograph of Duistermaat and Kolk \cite{DuKo}.
\end{example}

\begin{example}
Let us give an example of a non-integrable Lie algebroid.  Recall
(see, e.g., \cite{Mack1}) that any closed two-form $\omega\in
\Omega^2(M)$ defines a Lie algebroid structure on
$\A_{\omega}=TM\oplus\mathbb{L}$, where $\mathbb{L}=M\times\Rr$ is the
trivial line bundle over $M$, the anchor is $(X,\lambda)\mapsto X$ and
the Lie bracket is defined by
\[ [(X,f),(Y,g)]=([X, Y],X(g)-Y(f)+\omega(X,Y)).\]
Using the obvious splitting of $\A$, which has curvature
$\Omega_\sigma=\omega$, we see that the monodromy group based at $x$
is given by
\[ N_{x}(\A_{\omega})= \set{ \int_{\gamma} \omega:
     [\gamma]\in \pi_{2}(M, x)}\subset \mathbb{R} \]
and so coincides with the group of spherical periods of $\omega$. 
If this group is non-discrete we obtain a non-integrable Lie
algebroid. 

For example, on the 2-sphere $S^2$ denote by $\omega_{S^2}$ the
standard area form, and let $M=S^2\times S^2$ with the closed 2-form
$\omega_\lambda=\omega_{S^2}\oplus \lambda\omega_{S^2}$, where
$\lambda\in\Rr$. Then the monodromy group $N_{x}(\A_{\omega_\lambda})$
is discrete iff $\lambda$ is rational.

The reader will notice that in the symplectic case this
obstruction is Kostant's ``prequantization condition''.
\end{example}

\begin{example}
Let us consider the case of a Poisson manifold. In Poisson geometry,
$A$-paths are also called \emph{cotangent paths} (see
\cite{GiGo,Fer1}). Cotangent homotopies are given by an action of the
Lie algebra $P_{0}\Omega(M)$ of time-dependent 1-forms, vanishing at
the end-points, on the space of cotangent paths $P(T^*M)$.  The orbits
of this action have codimension $2\dim M$. In this case, this action
is the restriction of a Lie algebra action on the larger space of \emph{all}
paths $\tilde{P}(T^*M)=\set{a:I\to T^*M}$,  which is tangent to the
submanifold $P(T^*M)\subset\tilde{P}(T^*M)$. 

Since we have a natural identification $\tilde{P}(T^*M)\simeq
T^*P(M)$, where $P(M)$ denotes the Banach space of paths $\gamma:I\to
M$ which are piecewise of class $C^2$ (\footnote{We need cotangent
  paths to be piecewise of class $C^1$, so we require their base paths
  to be piecewise of class $C^2$. See also \cite{CrFe}.}), 
we have a natural symplectic structure on
$\tilde{P}(T^*M)$.  This brings symplectic geometry into the picture,
and we have the following result which is proved in \cite{CaFe}:

\begin{proposition}
The Lie algebra action of $P_{0}\Omega(M)$ on
$\tilde{P}(T^*M)$ is Hamiltonian with equivariant moment map
$J:\tilde{P}(T^*M)\to P_{0}\Omega(M)^*$ given by:
\[ \langle J(a),\eta\rangle=
\int_0^1 \langle  \frac{d}{dt}\pi(a(t))-\#a(t),\eta(t,\gamma(t))\rangle dt.\]
\end{proposition}

Since the level set $J^{-1}(0)$ is precisely the set of cotangent
paths $P(T^*M)$, it follows that in this case the Weinstein
groupoid can be described alternatively as a Marsden-Weinstein
reduction:
\begin{equation}
\label{eq:sigma:model}
 \G(T^*M)=\tilde{P}(T^*M)//P_{0}\Omega(M).
\end{equation}
The two alternative descriptions (\ref{eq:weinstein:groupoid}) and
(\ref{eq:sigma:model}) give the precise relationship between the
integrability approach introduced in \cite{CrFe}, which as we have
explained is valid for any Lie algebroid, and the approach of Cattaneo
and Felder in \cite{CaFe}, which is based on the Poisson sigma-model
and which only holds for Poisson manifolds.

Since we have the alternative description of the Weinstein groupoid of
a Poisson manifold as a Marsden-Weinstein reduction we obtain a
symplectic form on $\G(T^*M)$. If this groupoid is smooth this
symplectic form is compatible with the groupoid structure, and we
conclude that

\begin{theorem}[\textsc{Symplectic integration} \cite{CaFe,CrFe1}]
If the Weinstein groupoid $\G(T^*M)$ of a Poisson manifold $M$ is
smooth, then it is a symplectic groupoid integrating $M$.
\end{theorem}

Obviously, it is possible to form quotients of the groupoid $\G(T^*M)$
outside the symplectic category, to yield examples of groupoids
integrating $(M,\pi)$ and which are not symplectic
integrations. Results on the integrability of regular Poisson
manifolds are given in \cite{AlHe,DaHe}, and for a detailed discussion
we refer the reader to the forthcoming article \cite{CrFe1}.
\end{example}

\section{Holonomy for Lie algebroids}             %
\label{holonomy}                                  %
The theory of holonomy for algebroids relies, as in the case of
foliations, on the concept of \emph{connection}, since one wants to
compare the transverse Lie algebroid structure as we vary along a
leaf. Let us start by recalling the general notion of connection which
we have introduced in \cite{Fer1}.

\begin{definition}
Let $p:E\to M$ be a a fiber bundle over $M$, and
$\pi:A\to M$ a Lie algebroid over $M$. An \textsc{$A$-connection} 
on $E$ is a bundle map $h:p^*A\to TE$ which makes the
following diagram commute:
\[
\xymatrix{
p^*A\ar[r]^{h}\ar[d]_{\widehat{p}}& TE \ar[d]^{p_*} \\
A\ar[r]_{\#} &TM }
\]
\end{definition}

Note that in this definition $h:p^*A\to TE$ is a generalization of the
notion of \emph{horizontal lift} of tangent vectors that one finds in
the usual theory of connections: For any $a\in A$, $h(u,a)$ is the
horizontal lift of $a$ to the point $u\in E$ in the fiber over
$x=\pi(a)$, and the diagram means that $p_*h(u,a)=\#a$. Instead of
lifting tangent vectors in $TM$ we lift elements of $A$, the bundle
that replaces the tangent bundle.

Given some $A$-path $a:I\to A$, and a point $u_0\in E$ in the
fiber over the initial base point $x_0=\pi(a(0))$, we can
look for the horizontal lift $\gamma:I\to E$, i.~e., the
unique curve satisfying:
\[ \dot{\gamma}(t)=h(\gamma(t),a(t)), \qquad \gamma(0)=u_0.\]
Note that this system has a solution defined only for small $t\in
[0,\varepsilon)$. However, suppose that $E$ is a vector bundle and
that the connection has the property that $h(u,0)=0$. Then, it follows
from standard results in the theory of o.d.e.'s, that the solution
will be defined for all $t\in [0,1]$, provided we choose the initial
condition $u_0$ small enough.  In this way, we get a diffeomorphism
from a neighborhood of zero in the fiber $E_{x_0}$ over the initial
point onto a neighborhood of zero in the fiber $E_{x_1}$ over the final
point. Such a map is called, of course, \emph{parallel transport}
along the $A$-path $a:I\to A$.

Let us consider now a fixed leaf $i:L\hookrightarrow M$ of the Lie
algebroid $\pi:A\to M$.  We denote by $\nu(L)=T_LM/TL$ the normal
bundle to $L$ and by $p:\nu(L)\to L$ the natural projection.  By the
tubular neighborhood theorem, there exists a smooth immersion
$\tilde{i}:\nu(L)\to M$ satisfying the following properties:
\begin{enumerate}
\item[i)] $\tilde{i}|_Z=i$, where we identify the zero section $Z$ of
$\nu(L)$ with $L$;
\item[ii)] $\tilde{i}$ maps the fibers of $\nu(L)$ transversely to the
foliation of $M$;
\end{enumerate}
We shall define an $A_L$-connection on the normal bundle $\nu(L)$, so
that parallel transport for this connection will be the Lie algebroid
holonomy.

Assume that we have fixed such an immersion, and let $x\in L$. Each
fiber $F_x=p^{-1}(x)$ is a submanifold of $M$ transverse to the
foliation, and so we have the transverse Lie algebroid structure
$A_{F_x}\to F_x$ (see \cite{Fer1}).  Because $F_x$ is a linear space
we can choose a trivialization and identify the fibers $(A_{F_x})_u$
for different $u\in p^{-1}(x)$. Finally, we choose a complementary
vector subbundle $E\subset A$ to $A_{F_x}$:
\begin{equation}
\label{eq:decompose}
A_u=E_u\oplus (A_{F_x})_u.
\end{equation}
Note that, by construction, the anchor $\#:A\to TM$ maps $A_{F_x}$
onto $TF_x$, its restriction to $E$ is injective, and vectors in $\#E$
are tangent to the orbit foliation.

Let $\al\in A_x$. We decompose $\al$ according to
(\ref{eq:decompose}):
\[\al=\al^\parallel+\al^{\perp},\text{where } \al^\parallel\in E_x,
\quad \al^{\perp}\in (A_{F_x})_x.\]
For each $u\in F_x=p^{-1}(x)$, we denote by
$\tilde{\al}^\parallel_u\in E_u$ the unique element such that $d_u
p\cdot \#\tilde{\al}^\parallel_u=\#\al^\parallel$, and by
$\tilde{\al}^\perp_u\in (A_{F_x})_u$ the element corresponding to
$\al^{\perp}$ under the identification $(A_{F_x})_u\simeq
(A_{F_x})_x$. We also set $\tilde{\al}\equiv
\tilde{\al}^\parallel+\tilde{\al}^\perp$.

We can now define our connection: given $\al\in A_x$, $x\in L$,
and $u\in F_x$, the horizontal lift to $\nu(L)$ is the map
\[ h(u,\al)=\#\tilde{\al}_u\in T_u\nu(L).\]
By construction, we have the defining property of an $A$-connection:
\[ p_*h(u,\al)=\#\al,\quad u\in p^{-1}(x).\]
The definition of $h$ depends on several choices made: tubular
neighborhood, trivialization of $A_{F_x}$ and complementary vector
bundle $E$. The changes of choices will eventually lead to conjugate
holonomy homomorphisms (to be defined below).

Given an $A$-path $a:I\to A$ with initial point $x_0$ and final point
$x_1$, parallel transport gives us a diffeomorphism
$H_L(a)_0:F_{x_0}\to F_{x_1}$, defined on neighborhoods of the
origin. More is true: $H_L(a)_0$ is covered by a Lie algebroid
isomorphisms $H_L(a)$ from $A_{F_{x_0}}$ to $A_{F_{x_1}}$, so we have
\[
\xymatrix{
A_{F_{x_0}} \ar[r]^{H_L(a)}\ar[d]& A_{F_{x_1}}\ar[d] \\
F_{x_0}\ar[r]_{H_L(a)_0} &F_{x_1}}
\]
This can be seen as follows. Given any $A$-path $a(t)$ in $L$, we can
find a time-dependent section $\al_t$ of $A$ over $L$ such that
$\al_t(\gamma(t))=a(t)$, where $\gamma(t)=\pi(a(t))$. Now we can
define a time-dependent section $\tilde{\al}_t$ covering $\al_t$ such
that the horizontal lift of $a(t)$ is an integral curve of the
time-dependent vector field
\[ X_t=\#\tilde{\al_t},\]
so that $H_L(a)_0$ is the map induced by the time-1 flow of $X_t$ on
$F_{x_0}$.  Since the flow of $X_t$ is induced by the 1-parameter
family of Lie algebroid homomorphisms $\Phi_t^{\al_t}$ of $A$ obtained
by integrating the family $\tilde{\al}_t$, the homomorphisms
$\Phi_1^{\al_t}$ give a Lie algebroid isomorphism $H_L(a)$ from
$A_{F_{x_0}}$ to $A_{F_{x_1}}$, which covers $H_L(a)_0$.

Since $H_L(a)$ is the time-1 map of some flow it follows that if
$a'$ is another $A$-path in $L$ such that $x_1=x'_0$ we have
\begin{equation}
\label{eq:holonomy:homomorphism}
H_L(a \cdot a')= H_L(a)\circ H_L(a'),
\end{equation}
where the dot denotes concatenation of $A$-paths. We call $H_L(a)$
the \emph{$A$-holonomy} of the $A$-path $a(t)$. One extends the
definition of $H_L$ for piecewise smooth $A$-paths in the obvious way.

Denote by $\Agerm(A_{F_x})$ the group of germs at $0$ of Lie
algebroid automorphisms of $A_{F_x}$ which map $0$ to $0$, and
by $\Omega_A(L,x_0)$ the group of piecewise smooth
\emph{$A$-loops} based at $x_0$.

\begin{definition}
The \textsc{$A$-holonomy} of the leaf $L$ with base
point $x_0$ is the map
\[H_L:\Omega_A(L,x_0)\to\Agerm(A_{F_{x_0}}).\]
\end{definition}

Notice that the holonomy of a leaf $L$ depends on the tubular
neighborhood $\tilde{i}:\nu(L)\to M$, on the choice of trivialization,
and on the choice of complementary bundles.  However, two different
choices lead to conjugate homomorphisms.  This Lie algebroid holonomy
has, however, a major drawback: two $A$-paths whose base paths are
homotopic may have distinct holonomy. On the other hand, one can show,
that if they are homotopic as $A$-paths they lead to the same
holonomy, so we do get a homomorphism
\[H_L:\G(A)_x^x\to\Agerm(A_{F_{x}}).\]

For practical computations it is much more efficient to have a
homomorphism defined on the fundamental group $\pi_1(L,x)$.  In
\cite{Fer2}, following constructions given in \cite{Fer2} and
\cite{GiGo} for the Poisson case, we have introduced a notion of
\emph{reduced holonomy} which is homotopy invariant relative to the
base paths. Recall that $\Agerm(A_{F_x})$ denotes the group of germs
at $0$ of Lie algebroid automorphisms of $A_{F_x}$ which map $0$ to
$0$. By an \emph{inner Lie algebroid automorphism} of $A$ we mean an
automorphism which is the time-1 flow of some time-dependent
section.  We shall denote by $\Ogerm(A_{F_x})$ the corresponding group
of germs of \emph{outer Lie algebroid automorphisms} (\footnote{As usual, for a
Lie algebroid $A$, the group of outer Lie algebroid automorphisms is
the quotient $\Aut(A)/\Inn(A)$.}). 

In \cite{Fer2} we prove the following

\begin{proposition}
\label{prop:reduced:holonomy}
Let $x\in L\subset M$ be a leaf of $A$ with associated $A$-holonomy
$H_L:\Omega_A(L,x)\to\Agerm(F_{x})$. If $a_1(t)$ and $a_2(t)$ are
$A$-loops based at $x$ with base paths $\gamma_1\sim\gamma_2$
homotopic then $H_L(a_1)$ and $H_L(a_2)$ represent the same
equivalence class in $\Ogerm(F_{x})$.
\end{proposition}

Given a loop $\gamma$ in a leaf $L$ we shall denote by
$\bar{H}_L(\gamma)\in\Ogerm(A_{F_x})$ the equivalence class of
$H_L(a)$ for some piecewise smooth family $a(t)$ with
$\#a(t)=\gamma(t)$. The map
\[\bar{H}_L:\Omega(L,x)\to\Ogerm(A_{F_x})\]
will be called the \emph{reduced holonomy homomorphism} of $L$.  This
map extends to continuous loops and, by a standard argument, it
induces a group homomorphism $\bar{H}_L:\pi_1(L,x)\to\Ogerm(A_{F_x})$.

Recall that, for a foliation $\F$ of a manifold $M$, a \emph{saturated
set} is a set $S\subset M$ which is a union of leaves of $\F$. A leaf
$L$ is called \emph{stable} if it has arbitrarily small saturated
neighborhoods.  In the case of the orbit foliation of a Lie algebroid
a set is saturated iff it is invariant under all inner
automorphisms. Hence, a leaf is stable iff it is has arbitrarily small
neighborhoods which are invariant under all inner automorphisms.

We shall call a leaf $L$ \emph{transversely stable} if $N\cap L$ is a
stable leaf for the transverse Lie algebroid structure $A_N$, i.~e.,
if there are arbitrarily small neighborhoods of $N\cap L$ in $N$ which
are invariant under all inner automorphisms of $A_N$. Using this
notion of Lie algebroid holonomy one can prove a Reeb-type stability
theorem:

\begin{theorem}[\textsc{Stability Theorem} \cite{Fer2}]
\label{thm:local:stability}
Let $L$ be a compact, transversely stable leaf, with finite reduced
holonomy. Then $L$ is stable, i.e., $L$ has arbitrarily small
neighborhoods which are invariant under all inner
automorphisms. Moreover, each leaf near $L$ is a bundle over $L$ whose
fiber is a finite union of leaves of the transverse Lie algebroid
structure.
\end{theorem}

The local splitting theorem for Lie algebroids (see \cite{Fer2},
Theorem 1.1) states that locally a Lie algebroid splits as a product
of $A_L$ and $A_N$.  Using holonomy one can investigate if a
neighborhood of a leaf trivializes as in the local splitting
theorem. An obvious necessary condition is that the Lie algebroid
holonomy be trivial. We refer the reader to \cite{Fer1,Fer2} for
further details.

\section{Linear Holonomy}                         %
\label{linear:holonomy}                           %

Recall (see \cite{Fer2,Wein1}) that a \emph{linear Lie algebroid} is
a Lie algebroid $\pi:A\to V$ such that:
\begin{enumerate}
\item[(i)] The base is a vector space $V$ (so $\pi:A\to V$ is
  trivial);
\item[(ii)] For any trivialization, the bracket of
  constant sections is a constant section;
\item[(iii)] For any trivialization and constant section $\al$, the
  vector field $\#\al$ is linear.
\end{enumerate}
The reader should notice that a linear Lie algebroid $A$ is isomorphic
to a (linear) action Lie algebroid $\gg\times V$.

For a Lie algebroid $A$, the transverse Lie algebroid structure to a
leaf of $A$ is a germ of a Lie algebroid for which the anchor vanishes 
at the base point. One should expect then that there is a well-defined
linear approximation. In fact, we have the following proposition:

\begin{proposition}
Let $A\to M$ be any Lie algebroid, fix $x_0\in M$ and denote by $L$
the leaf through $x_0$. There is a
natural linear Lie algebroid structure $A^{\text{lin}}$ over 
the normal space $N_{x_0}=T_{x_0}M/T_{x_0}L$ with
\[ A^{\text{lin}}=\gg_{x_0}\times N_{x_0}.\]
\end{proposition}

We call $A^{\text{lin}}\to N_{x_0}$ the \emph{linear approximation} to $A$ at
$x_0$. We also have a linear version of holonomy 
\[H_L^{\text{lin}}\equiv dH_L:\Omega_A(L,x)\to 
\Aut(\gg)\times \textrm{GL}(F_x)\simeq \Aut(A^{\text{lin}}).\]
Moreover, we can obtain linear holonomy as parallel transport along a linear
$A$-connection generalizing the Bott connection of ordinary
foliation theory, and which we will now recall.

First, if $p:E\to M$ is a vector bundle over $M$, a linear $A$-connection on
$E$ is a connection $h:p^*A\to TE$ which is linear:
\[ h(u,\al_1+\al_2)=h(u,\al_1)+h(u,\al_2),\qquad
h(u,\lambda\al)=\lambda h(u,\al),\]
where $\lambda\in\Rr$, $u\in E$ and $\alpha_i\in A$.

In exactly the same way as one does with ordinary connections, we can
associate with a linear connection an operator 
$\nabla:\Gamma(A)\times\Gamma(E)\to\Gamma(E)$ which satisfies 
\begin{enumerate}
\item[(i)] $\nabla_{\al_1+\al_2}s=\nabla_{\al_1} s+\nabla_{\al_2} s$;
\item[(ii)] $\nabla_\al(s_1+s_2)=\nabla_\al s_1+\nabla_\al s_2$;
\item[(iii)] $\nabla_{f\al}s=f\nabla_\al s$;
\item[(iv)] $\nabla_\al(f s)=f\nabla_\al s+\#\al(f)s$;
\end{enumerate}
where $\al,\al_1,\al_2\in\Gamma(A)$, $s,s_1,s_2\in\Gamma(E)$,
and $f\in C^{\infty}(M)$. Conversely, any operator satisfying (i)-(iv)
determines a linear connection. For a detailed discussion of such
connections see \cite{Fer2}.

The \emph{Bott connection} of a Lie algebroid is a pair of linear
$A$-connections on $\Ker\#|_L$ and on $\nu^*(L)$. On one hand, we have 
on the vector bundle $\Ker\#|_L$ the linear $A$-connection: 
\begin{equation}
\label{Bott:connection}
\nabla^L_\al\gamma\equiv {[\tilde{\al},\tilde{\gamma}]},
\end{equation}
where $\tilde{\al},\tilde{\gamma}\in\Gamma(A)$ are sections extending
the sections $\al\in \Gamma(A_L)$ and $\gamma\in\Gamma(\Ker\#|_L)$. On
the other hand, we have on the conormal bundle $\nu^*(L)=\set{\omega\in
  T^*_LM:\omega|TL=0}$ the linear $A$-connection:
\begin{equation}
\label{Bott:connection:*}
\check{\nabla}^L_\al\omega\equiv
\Lie_{\#\tilde{\al}}\tilde{\omega}|_x,
\end{equation}
where now we take a section $\tilde{\al}\in\Gamma(A)$ and a 1-form
$\tilde{\omega}\in \Omega^1(M)$ extending the sections $\al\in
\Gamma(A_L)$ and $\omega\in\Gamma(\nu^*(L))$. It is easy to check that
$\nabla^L$ and $\check{\nabla}^L$ satisfy properties (i)-(iv) above. 

It is more convenient to consider the connections $\nabla^L$ and
$\check{\nabla}^L$ together, rather than leaf by leaf. Also, we
can work on the direct sum $A\ominus TM$ (the reason for the strange
symbol will be explained later) and replace the pair of connections by
a single connection, so we set:

\begin{definition}
\label{defn:basic:connection}
A linear connection $\nabla=\nabla^A\ominus\nabla^M$ on $A\ominus TM$ 
is called a \textsc{basic connection} if
\begin{enumerate}
\item[i)] $\nabla$ is compatible with the Lie algebroid structure, i.e.,
\[\nabla^M\#=\#\nabla^A;\]
\item[ii)] $\nabla^A$ and $\nabla^M$ restrict to the Bott connection on each
  leaf $L$, i.~e., if $\al,\gamma\in\Gamma(A)$, $\omega\in\Omega^1(M)$,
  with $\#\gamma|_L=0$ and $\omega|_{TL}=0$, then
\[ \nabla_\al(\gamma,\omega)|_L=([\al,\gamma],\Lie_{\#\al}\omega)|_L.\]
\end{enumerate}
\end{definition}

Basic connections always exist. A simple procedure, due to Crainic, for
constructing a basic connection is to start with any $TM$-connection
$\tilde{\nabla}$ on $A$ and set:
\[ \nabla_\al(\gamma,X)=(\tilde{\nabla}_{\#\gamma}\al+[\al,\gamma],
\#\tilde{\nabla}_X\al+[\#\al,X]).\]

The holonomy along a leaf $L$ of a basic connection $\nabla$ gives the
linear holonomy of $L$ introduced above in the following way: the
holonomy of the basic connection $\nabla$ determines endomorphisms of
the fiber $A_x$ which map $\ker\#_x$ isomorphically into itself, and
these are the linear holonomy maps. Moreover, we have the following
Bott-type vanishing theorem which lies at the basis of the secondary
characteristic classes to be introduced in the next section.

\begin{theorem}[\textsc{Bott vanishing theorem} \cite{Fer1}]
Let $R$ denote the curvature of a basic connection. Then for any
section $\gamma\in\Gamma(A)$ and 1-form $\omega\in\Omega^1(M)$
satisfying $\#\gamma|_L=0$ and $\omega|_{TL}=0$, we have
\[ R(\al,\gamma)(\gamma,\omega)|_L=0.\]
\end{theorem}

For more on Bott vanishing theorems in the special context of regular Lie
algebroids we refer to the work of Kubarski (\cite{Ku}).

\section{Characteristic classes}                  %
\label{ch:classes}                                %

One can define (primary) characteristic classes for algebroids as one
does in the usual Chern-Weil theory. For example, let $E\to M$ be a
real vector bundle with rank $q$ and pick some $A$-connection $\nabla$
on $E$. Its curvature
\[
R(\al,\beta)=\nabla_\al\nabla_\beta-\nabla_\beta\nabla_\al+
\nabla_{[\al,\beta]},
\] 
defines a linear map $R_{\al,\beta}=R(\al,\beta):E_x\to E_x$ which
satisfies $R_{\al,\beta}=-R_{\beta,\al}$. Hence, the map
$(\al,\beta)\mapsto R_{\al,\beta}$ can be considered as a
$\gl(E)$-valued 2-section, and by fixing a basis of local sections for
$E$, so that $E_x\simeq\Rr^q$, we have that
$R_{\al,\beta}\in\gl_q(\Rr)$. This matrix representation of
$R_{\al,\beta}$ is defined only up to a change of basis in
$\Rr^q$. Therefore, if
\[ P:\gl_q(\Rr)\times\cdots\times\gl_q(\Rr)\to\Rr\]
is a symmetric, $k$-multilinear, $\Ad(GL_q(\Rr))$-invariant function,
we can introduce a well-defined $2k$-form
$\lambda(\nabla)(P)\in\Omega^{2k}(A)$ by the formula
\begin{equation}
\label{eq:Chern-Weil:homomorphism:connection}
\lambda(\nabla)(P)(\al_1,\dots,\al_{2k})=\sum_{\sigma\in S_{2k}}
(-1)^{\sigma}
P(R_{\al_{\sigma(1),\sigma(2)}},\dots,R_{\al_{\sigma(2k-1),\sigma(2k)}}).
\end{equation}
This form is closed and hence defines a certain Lie algebroid cohomology class
$[\lambda(\nabla)(P)]\in H^{2k}(A)$. It is not hard to see that this
cohomology class is independent of the choice of connection, so we
have defined some intrinsic characteristic classes of the vector
bundle $E$. For example, if we let $P_k$ be the elementary symmetric
polynomials we have the $A$-Pontrjagin classes
\[ p_k(E,A)=[\lambda(P_{2k})]\in H^{4k}(A).\]
As usual, one does not need to consider the classes for odd $k$ since
we have
\[ [\lambda(P_{2k-1})]=0,\]
as can be seen by choosing a connection compatible with a Riemannian
metric.

However, these classes do not really contain any new information. In
fact, the anchor $\#:A\to TM$ determines a chain map 
$\#^*:(\Omega^\bullet(M),d)\to(\Omega^\bullet(A),d_A)$ and so we have
an induced map in cohomology:
\[ \#^*:H^\bullet(M)\to H^\bullet(A).\]
If we choose some ordinary connection $\tilde{\nabla}$ on $E$ and take
$\nabla_\al=\tilde{\nabla}_{\#\al}$ we see immediately that
\[ p_k(E,A)=\#^*p_k(E),\]
where $p_k(E)$ are the usual Pontrjagin classes of $E$. The same is
true for any (primary) characteristic classes one may define: we have
a commutative diagram
\[
\xymatrix{
I^\bullet(G)\ar[r]\ar[dr]& H^\bullet_{{\text{de Rham}}}(M) \ar[d]^{\#^*} \\
& H^\bullet(A)}
\]
where on the top row we have the usual Chern-Weil homomorphism and
on the diagonal we have the $A$-Chern-Weil homomorphism (see
\cite{Fer1}). These classes were introduced for Poisson manifolds by
Vaisman in \cite{Vais1}, and for regular Lie algebroids by Kubarski in
\cite{Ku2}.

The fact that all these classes arise as image by $\#^*$ of well-known
classes is perhaps a bit disappointing. However, one can define
\emph{secondary characteristic classes} which are true invariants of the Lie
algebroid, in the sense that they do not arise as images by $\#^*$ of
some de Rham cohomology classes. These classes are analogous to
the exotic classes of foliation theory introduced by Bott \emph{et
al.}~(see, e.g.,~\cite{Bott}). To define them introduce a pair of
connections $(\nabla^1,\nabla^0)$ on $A\ominus TM$ where:
\begin{itemize}
\item $\nabla^1$ is a basic connection;
\item $\nabla^0$ is a Riemannian connection (i.e., 
$\nabla^0_\al=\nabla_{\#\al}$ with $\nabla$ the Levi-Civita
connection);
\end{itemize}
Given an $\Ad$-invariant, symmetric polynomial $P$, the classes
are defined by a transgression formula in the spirit of Chern and Simons
(\cite{Chern}): 
\begin{multline}
\label{eq:invariants:2}
\lambda^{1,0}(P)(\al_1,\dots,\al_{2k-1})=\\
k \sum_{\sigma\in S_{2k-1}} (-1)^{\sigma}\int_{0}^1
P(\nabla^{1,0}_{\al_{\sigma(1)}},R^t_{\al_{\sigma(2)},\al_{\sigma(3)}},
\dots,R^t_{\al_{\sigma(2k-2)},\al_{\sigma(2k-1)}})dt,
\end{multline}
where $\nabla^{1,0}=\nabla^1-\nabla^0$ and $R^t$ is the curvature
of $\nabla^t=(1-t)\nabla^0-t\nabla^1$. Again we have:
\begin{theorem}[\textsc{Secondary classes} \cite{Fer1}]
Let $k$ be odd. Then
\begin{enumerate}
\item[(i)] $\lambda^{1,0}(P)\in\Omega^{2k-1}(A)$ is closed;
\item[(ii)] The cohomology class $[\lambda^{1,0}(P)]\in H^{2k-1}(A)$
  is independent of the choice of connections;
\end{enumerate}
\end{theorem}

In general these classes do not lie in the image of
$\#^*:H^\bullet(M)\to H^\bullet(A)$, as can be seen from some of the examples
given below, so we obtain genuine invariants of the Lie
algebroid. In particular, if we take $P=P_k$ the elementary symmetric
polynomials we obtain the \emph{secondary characteristic classes} of a
Lie algebroid:
\[ m_k(A)=[\lambda^{1,0}(P_k)]\in H^{2k-1}(A), \qquad k=1,3,5,\dots\]
Explicit computation of these classes appear in \cite{Fer1}. The best
understood class is $m_1$ and it was known before \cite{Fer1} as the
\emph{modular class} of a Lie algebroid. This class was introduced
first by Weinstein in \cite{Wein2}, for the special case of Poisson
manifolds, and by Weinstein \emph{et al.}~in \cite{ELW} for general
Lie algebroids. There the following geometric interpretation was
given. Let us think of sections of the line bundle
$Q_A=\wedge^{\text{top}}A\otimes \wedge^{\text{top}}TM$ (or
$Q_A\otimes Q_A$ in the non-orientable case) as ``transverse
measures'' in $A$. Then the modular class is the obstruction to the
existence of invariant transverse measures: $m_1(A)=0$ iff there
exists a measure invariant under the flow of any section of $A$. In
the Poisson case, $m_1(A)=0$ iff there exists a (true) measure
invariant under the flow of any hamiltonian diffeormorphism. The
modular class was also studied in a purely algebraic context by
Huebschmann \cite{Hu1,Hu2} and Xu \cite{Xu1}.

Crainic in \cite{Cra} has developed, independently from \cite{Fer1},
a theory of characteristic classes of representations. Recall that a
\emph{representation} of a Lie algebroid $A\to M$ is a vector bundle $E\to M$
together with a map $\Gamma(A)\times\Gamma(E)\to\Gamma(E)$ satisfying:
\begin{enumerate}
\item[(i)] $(f\al)\cdot s=f\al\cdot s$;
\item[(ii)] $\al\cdot (fs)=f\al\cdot s+\#\al(f)s$;
\item[(iii)] $[\al,\beta]\cdot s=\al\cdot(\beta \cdot
  s)-\beta\cdot(\al \cdot s)$; 
\end{enumerate}
If we set $\nabla_\al s\equiv\al\cdot s$ then we see that a
representation is nothing other than a flat $A$-connection on $E$. The
terminology is motivated by the case when $A$ is a Lie algebra.

Assume then that $E\to M$ is a representation of a Lie algebroid
$A$ and denote by $\nabla^1$ the associated flat $A$-connection on
$E$. Then $\nabla^1$ induces an adjoint connection $\nabla^{1*}$ on 
the dual bundle $E^*$. If we pick some Riemannian metric on $E$ we
obtain an identification $E\simeq E^*$, so that the adjoint connection
determines a new $A$-connection $\nabla^0$ on $E$. We
can then check that if $P$ is some $\Ad$-invariant polynomial,
the transgression formula (\ref{eq:invariants:2}) gives characteristic
classes of the representation(\footnote{The classes are actually
  defined only for complex representations and for polynomials
  representing classes in the relative cohomology
  $H^\bullet(GL(n),U(n))$, but we shall ignore this aspect here.}):
\[ u(E,P)=[\lambda^{1,0}_E(P)]\in H^\bullet(A).\]
These are the characteristic classes introduced by Crainic in
\cite{Cra}. Their main properties are given in the following
proposition:

\begin{proposition}
For representations $E$ and $F$ of a Lie algebroid $A$, and any
characteristic class $u(\cdot)=u(\cdot,P)$ as defined above, we have:
\begin{enumerate}
\item[(i)] $u(E\oplus F)=u(E)+u(F)$;
\item[(ii)] $u(E\otimes F)=\rank(E)u(F)+\rank(F)u(E)$;
\item[(iii)] $u(E^*)=-u(E)$;
\end{enumerate}
\end{proposition}

In particular, the characteristic classes vanish if $E$
admits an invariant metric, hence these classes measure the obstruction
to the existence of such invariant metrics.

Notice that the Crainic classes are representation dependent, while
the classes introduced above were intrinsic classes of the Lie
algebroid. If there existed some natural or canonical representation
of a Lie algebroid $A$, one would expect that the (extrinsic) classes
of such representation would give the intrinsic characteristic classes
of $A$. However, a general Lie algebroid has no adjoint action as
opposed to, say, a Lie algebra. It turns out that it is still possible
to recover the intrinsic characteristic classes from characteristic
classes of a representation if one weakens the notion of
representation to a \emph{representation up to homotopy}. Then there
is a natural \emph{adjoint representation up to homotopy} for every
Lie algebroid $A$, and the Crainic classes of this representation
yield the intrinsic characteristic classes.

The notion of connection up to homotopy is obtained as a weaker version
of Quillen's (see \cite{Qu}) notion of superconnection. Recall that a
super-vector bundle is just a $\Zz_2$-graded vector bundle over a
manifold $M$. If $(E,\partial)$ is a super-complex of vector
bundles over $M$:
\[ \xymatrix{
  E^1\ar@<-.5ex>[r]_{\partial_1}&E^0\ar@<-.5ex>[l]_{\partial_0}} \]
we can view it as an element in the $K$-theory of $M$,
i.e., the formal differences $E=E^0\ominus E^1$. We set:
\begin{definition}
An \textsc{$A$-connection up to homotopy} on a super-vector bundle
$(E,\partial)$ is a $\Rr$-bilinear map
$\nabla:\Gamma(A)\times\Gamma(E)\to \Gamma(E)$, such that:
\begin{enumerate}
\item[(i)] $\nabla$ perserves the grading and
$\nabla\partial_i=\partial_{i+1}\nabla$;
\item[(ii)] $\nabla_\al(f s)=f\nabla_\al s+\#\al(f)s$;
\item[(iii)] $\nabla_{f\al}s=f\nabla_\al s+[H(f,\al),\partial]s$;
\end{enumerate}
where the homotopy $H(f,\al):E\to E$ has degree 1.
\end{definition}

A flat connection up to homotopy is the same as a \emph{representation
up to homotopy} (see \cite{ELW}, for details and background on these
connections). For representations up to homotopy one can define extrinsic
characteristic classes $[\lambda^{1,0}_E(P)]\in H^\bullet(A)$ in a similar
manner to the case of ordinary representations (see \cite{Cra2}).
The main example we are interested in is the \emph{adjoint representation
up to homotopy} of a Lie algebroid $A$. It is the representation
\[ \ad(A):  \xymatrix{ A\ar@<-.5ex>[r]_{\#}&TM\ar@<-.5ex>[l]_{0}}\]
where we take
\[ \nabla^{\ad(A)}_\al(\gamma,X)=([\al,\gamma],\Lie_{\#\al}X).\]
The reader may check that this is indeed a flat connection up to
homotopy with homotopy map $H(f,\al)$ given by:
\[ H(f,\al)\gamma=0, \qquad H(f,\al)X=X(f)\al.\]

At this point the reader will notice the close relationship between
the adjoint representations up to homotopy and the basic connections
introduced in the previous section: every basic connection determines
a flat connection up to homotopy on the super-vector bundle $A\ominus
TM$. To make this relationship more precise, and finally justify the
use of the symbol $\ominus$ as a ``difference'' similar to K-theory,
we introduce a notion of equivalence among connections up to homotopy:

\begin{definition}
\label{def:equivalence:rep}
Two connections up to homotopy $\nabla$ and $\nabla'$ on a
super-vector bundle $E$ are said to be \textsc{equivalent} if there
exists some degree-zero $\End(E)$-valued 1-form $\theta$ such that:
\[ \nabla'_\al=\nabla_\al+[\theta(\al),\partial].\]
\end{definition}

One can then show (see \cite{Cra2} for details):

\begin{proposition}
A true $A$-connection $\nabla$ on $A\ominus TM$ is equivalent to
$\nabla^{\ad(A)}$ iff $\nabla$ is a basic connection.
\end{proposition}

Finally, Crainic has established that the intrinsic
characteristic classes of a Lie algebroid coincide with the
characteristic classes of the adjoint representation:
\[ [\lambda^{1,0}_{\ad(A)}(P)]=[\lambda^{1,0}(P)].\]
Namely, he shows that, in general, two equivalent representations up
to homotopy have the same characteristic classes.  We refer the reader
to \cite{Cra2} for proofs and further details.

\section{K-theory}                                %
\label{ktheory}                                   %
It is well-known that $K$-theory is the most efficient of all
cohomology theories admitting a geometric description. Since we have a
notion of representation of a Lie algebroid it is therefore natural to
look at the possibility of constructing a $K$-theory in the context of
Lie algebroids. A first attempt in constructing a $K^0$-functor was
made by Ginzburg in \cite{Ginz} and we recall his construction in
this section. This should by no means be considered as the final word
on this theory. In fact, we have good indications that a $K$-theory
based on the notion of representations up to homotopy would be more
efficient, and this is the subject of current investigation.

We denote by $\Vect_A(M)$ the semi-ring of equivalence classes of
representations of a Lie algebroid $\pi:A\to M$. The following example shows
that this semi-ring is usually too large.

\begin{example}
\label{ex:trivial:alagebroid}
Consider the trivial Lie algebroid $\pi:A\to M$ of rank $r$, and let
$p:E\to M$ be an ordinary vector bundle over $M$. Any family
$\set{\sigma_x}_{x\in M}$ of representations of the
abelian Lie algebras $A_x$ on the vector space $E_x$,
\[ \sigma_x:A_x\to\End(E_x),\]
determines a representation of $A$:
\[ \nabla_\al s|_x=\sigma_x(\al(x))\cdot s(x).\]
\end{example}

We must therefore somehow reduce the number of relevant
representations. Ginzburg in \cite{Ginz} suggests reducing the number
by identifying representations which can be deformed into one another.

\begin{definition}
Let $E_0$ and $E_1$ be representations of a Lie algebroid $A$. We say
that $E_0$ and $E_1$ are \textsc{deformation equivalent} if they are
isomorphic to representations that can be connected by a family $E_t$
of representations of $A$.
\end{definition}

This clearly defines an equivalence relation on $\Vect_A(M)$ which
respects the semi-rimg structure. Denote then by $\Vectd_A(M)$ the
semi-ring of equivalence classes of representations of a Lie algebroid
$A$.

\begin{example}
For a trivial Lie algebroid $\pi:A\to M$, all
representations are homotopic to one another so
$\Vectd_A(M)=\Vect(M)$. In fact, if $p:E\to M$ is some representation
with an associated flat connection $\nabla$, then $\bar{\nabla}^t\equiv
t\nabla$ defines a family of flat connections giving a deformation of
$(E,\nabla)$ to the trivial representation.
\end{example}

Although this is not the only possible way of reducing equivalence
classes of representations (see e.g.~\cite{Ginz} for the related concept of
\emph{homogeneous} representations), we restrict our attention to the
deformation equivalent classes, for this already gives us the flavor
of any such theory.

\begin{definition}
The \textsc{$K$-ring} of $A$ is the Grothendieck ring $K(A)$ associated
with the semi-ring $\Vect_A$.
\end{definition}

To check that these rings are reasonable objects let us look at some examples:

\begin{example}
In example \ref{ex:trivial:alagebroid}, we saw that for a trivial Lie
algebroid $A$ the $K$-ring $K(A)$ is equal to the ordinary $K$-ring of
$M$, i.e., $K(A)=K(M)$.
\end{example}

\begin{example} 
Let $A=TM$. Then a representation $E\to M$ of $A$ is just a vector
bundle together with an (ordinary) flat connection $\nabla$, or
equivalently, a representation of $\pi_1(M)$. Two
representations are deformation equivalent iff the corresponding
representations of $\pi_1(M)$ lie in the same
path-connected component of the space of representations of
$\pi_1(M)$. Hence, we conclude that $K(TM)=\pi_0(\text{Rep}(\pi_1(M)))$.

\end{example}

\begin{example}
Let $A=\gg^*\times\gg\to\gg^*$ be the transformation Lie algebroid 
associated with the coadjoint action
$\ad^*:\gg\to\mathfrak{gl}(\gg^*)$. Also, let $p:E\to \gg^*$ be some
representation of $A$, with an associated flat connection $\nabla$.

There is a Lie algebra representation $\rho:\gg\to V$, where
$V=p^{-1}(0)$, naturally associated with the representation. It can be
defined as follows: if $y\in\gg$ and $v\in V$ choose sections
$\al\in\Gamma(A)$ and $s\in\Gamma(E)$ such that $\al(0)=y$ and
$s(0)=v$. Then
\[ \rho(y)\cdot v\equiv \nabla_\al s(0),\]
and it is easy to check that this is independent of all choices. 

The Lie algebra representation $\rho$ defines a new representation
$(\bar{E},\bar{\nabla})$ of the Lie algebroid $A$: $\bar{E}$ is the
trivial vector bundle $\gg^*\times V\to \gg^*$, and $\bar{\nabla}$ is
the unique flat connection which for a constant section $s(x)=v$
satisfies
\[ \bar{\nabla}_\al s=\rho(\al)\cdot v,\]
(here we identify $\al\in\Gamma(A)$ with a function $\al:\gg^*\to\gg$).

One can check (see \cite{Ginz}) that the representations $(E,\nabla)$ and
$(\bar{E},\bar{\nabla})$ are deformation equivalent. Moreover, two
representations $(E_0,\nabla^0)$ and $(E_1,\nabla^1)$ are deformation
equivalent iff the associated Lie algebra representations $\rho_0$ and
$\rho_1$ are in the same path component of the space of
representations of $\gg$. Hence, we conclude that $K(\gg^*)=R(\gg)$,
the ring of representations of $\gg$.
\end{example}

To complete the properties of this $K$-theory we consider \emph{Morita
equivalence} of Lie algebroids. The definition is based on the notion
of \emph{pull-back Lie algebroid} due to Higgins and Mackenzie (see
\cite{HiMa}): we start with a Lie algebroid $A\to M$ and we consider a
surjective submersion $\phi:Q\to M$. Then the pull-back Lie algebroid
$\phi^\star A$ completes the diagram
\[
\xymatrix{ \phi^\star A\ar@{-->}[r]^{\widehat{\phi}}\ar@{-->}[d]& A\ar[d]\\
Q \ar[r]^\phi & M}
\]
so that $\widehat{\phi}:\phi^\star A\to A$ is a morphism of Lie algebroids.
As a vector bundle, $\phi^\star A$ is given by:
\[ \phi^\star A=\set{(\al,X)\in \phi^*A\times TQ:\#\al=d\rho\cdot
  X},\]
where $\phi^*A$ is the usual pull-back of vector bundles. For the
anchor one takes projection into the second factor, while the bracket
is given by:
\[ [(f\al,X),(g\beta,Y)]=(fg[\al,\beta]+X(g)\al-Y(f),[X,Y]),\]
whenever $\al,\beta$ are pull-backs of sections of $A$, $f,g\in
C^\infty(Q)$ and $X,Y\in\X(Q)$. 
In the sequel we shall use $\phi^\star A$ to denote the pull-back
construction in the category of Lie algebroids and $\phi^*A$ to denote
the pull-back construction in the category of vector bundles. Since we
have $\phi^\star A\simeq\phi^* A\oplus \Ker\phi_*$, these constructions
only coincide if $\phi$ is a covering. 

\begin{definition}
Two Lie algebroids $A_1\to M_1$ and $A_2\to M_2$ are called \textsc{Morita
equivalent} if there exists a pair of surjective submersions, with
simply connected fibers,
\[ \xymatrix{
  &Q\ar[dl]_{\phi_1}\ar[dr]^{\phi_2}\\M_1& &M_2} \]
such that the pull-back Lie algebroids $\phi_1^\star A_1$ and
$\phi_2^\star A_2$ are isomorphic.
\end{definition}

\begin{remark}
Morita equivalence was first introduced in the context of Poisson
manifolds by Xu in \cite{Xu2}, and further studied by Ginzburg and Lu
in \cite{GL}. This definition makes sense only for \emph{integrable}
Poisson manifolds. For Lie algebroids, the definition above is due to
Ginzburg \cite{Ginz}, and is a linear version of the notion of Morita
equivalence of Lie groupoids (see e.g.~\cite{Cra}).  For Poisson
manifolds, this notion of Morita equivalence (applied to the cotagent
bundle Lie algebroids) is weaker than Xu's definition of Morita
equivalence, but it makes sense for \emph{all} Poisson manifolds. For
this reason it is called in \cite{Ginz} \emph{weak Morita
equivalence}. The advantage of this definition is that the invariants
we have been discussing are in fact (weak) Morita invariants. The
notion of weak Morita equivalence is just one of several
possibilities (see the discussion in \cite{Ginz}).
\end{remark}

Let $E$ be a representation of the Lie algebroid $A$. Then $\phi^*E$
is naturally a representation of $\phi^\star A$: the flat connections
on $E$ and $\phi^*E$ are related by
\[ \tilde{\nabla}_{(\al,X)}(fs)=X(f)s+f\nabla_\al s,\]
whenever $s$ is the pull-back of a section of $E$. Conversely, since
the $\Ker d\phi$-action on $\phi^*E$ coincides with the action
obtained from the natural flat connection on $\phi^*E$ along the 
$\phi$-fibers, it follows that pull-back induces a bijection
\[ \Vect_A(M)\longleftrightarrow \Vect_{\phi^\star A}(Q).\]
Therefore, the definition implies that Morita equivalent Lie
algebroids have the same representations:

\begin{theorem}[\textsc{Invariance under Morita Equivalence} \cite{Cra,Ginz}]
Let $\pi_1:A_1\to M_1$ and $\pi_2:A_2\to M_2$ be Morita equivalent Lie
algebrois. Then $\Vect_{A_1}(M_1)\simeq \Vect_{A_2}(M_2)$ and 
$\Vectd_{A_1}(M_1)\simeq \Vectd_{A_2}(M_2)$. In particular, we have:
\[ K(A_1)\simeq K(A_2).\]
\end{theorem}

One can also show, that Morita equivalent Lie algebroids have
isomorphic zero and first Lie algebroid cohomology groups (see
\cite{Ginz}). If one assumes further that the $\phi_i$-fibers are
$n$-connected, then one can show that they have isomorphic Lie
algebroid cohomologies $H^k(A_1)\simeq H^k(A_2)$, for all $k\le n$
(see \cite{Cra,Ginz}).

Most of, if not all, the properties of Morita equivalence, can be
reduced to the statement that Morita equivalent Lie algebroids have
the same orbit space. In fact, the algebroids $\phi^\star
A_1$ and $\phi^\star A_2$ obviously have isomorphic foliations and
isomorphic transverse structures. Hence, the assignment $L\mapsto
\phi_2(\phi^{-1}_1(L))$ gives a bijection between the leaves of $A_1$
and the leaves of $A_2$. Moreover, if $x_1=\phi(q)$ and $x_2=\phi(q)$, 
it follows from the pull-back construction that the transverse Lie
algebroid structures to $A_1$ at $x_1$ and to $A_2$ at $x_2$, are 
isomorphic. Hence:

\begin{theorem}
Let $A_1$ and $A_2$ be Morita equivalent Lie algebroids. 
Then there exists a 1:1 correspondence between the leaves of $A_1$ and
the leaves of $A_2$. Moreover, the transverse Lie algebroid
structures of corresponding leaves are isomorphic. 
\end{theorem}

At this stage it is not clear how one can extend the results of this
section to representations up to homotopy, an aspect of Lie algebroid
theory we feel deserves further investigation.

\bibliographystyle{amsplain}
\def\lllll{}

\end{document}